\documentclass[12pt]{amsart}

\title{AF-embeddings of residually finite dimensional C*-algebras}
\title{AF-embeddings of residually finite dimensional C*-algebras}
\author{Marius Dadarlat}\address{MD: Department of Mathematics, Purdue University, West Lafayette, IN 47907, USA}\email{mdd@purdue.edu}	
\thanks{M.D. was partially supported by NSF grant \#DMS--1700086}
\urladdr{http://www.math.purdue.edu/~mdd/}
\usepackage{enumerate}
\usepackage{amssymb}
\newcommand{\ep}{\varepsilon}

\newcommand{\Z}{\mathbb{Z}}
\newcommand{\Q}{\mathbb{Q}}

\theoremstyle{plain}
\newtheorem{theorem}[subsection]{Theorem}

\newtheorem{proposition}[subsection]{Proposition}
\theoremstyle{definition}
\newtheorem{remark}[subsection]{Remark}
\newtheorem{definition}[subsection]{Definition}

\begin{document}

\begin{abstract}
It is shown that a separable exact residually finite dimensional C*-algebra with locally finitely generated  (rational) $K^0$-homology embeds in a uniformly hyperfinite $C^*$-algebra.
\end{abstract}

\maketitle

\section{Introduction}

Kirchberg  proved that any separable exact C*-algebra embeds in the Cuntz algebra $\mathcal{O}_2$, see \cite{Kir-icm}. A related major open problem asks if any separable exact and quasidiagonal C*-algebra
embeds in an almost finite dimensional algebra (AF-algebra), see \cite[Ch.8]{Brown-Ozawa}.
Most positive results on AF-embeddability depend on the universal coefficient theorem in KK-theory (abbreviated UCT) \cite{RosSho:UCT}, see for example \cite{Dad:kk-top}, \cite{ORS:elem-amen}, \cite{TWW:quasidiagonality}.
A  general result of Ozawa \cite{Ozawa:qd} shows that the cone over an exact separable C*-algebra  is AF-embeddable.
Such cones are automatically quasidiagonal by a theorem of Voiculescu \cite{Voi:qd}.
While Ozawa's proof does not use the UCT explicitly, cones are contractible and in particular they do satisfy the UCT. Cones also play a key role in R{\o}rdam's paper on purely infinite AH-algebras and AF-embeddings \cite{Rordam:AF}. Indeed, R{\o}rdam's $C^*$-algebra $\mathcal{A}[0,1]$, which he showed that contains the cone over $\mathcal{O}_2$ as a subalgebra, is itself an inductive limit of cones over matrix algebras and in particular it is KK-contractible.
In a very recent paper  \cite{Jamie-Gabe:connective}, Gabe proves  that a separable exact $C^*$-algebra for which its primitive spectrum has no non-empty compact open subsets embeds in $\mathcal{A}[0,1]$ and hence it is AF-embeddable. As far as I am aware, all previously known AF-embeddings results which do not assume the UCT factor through  inductive limits of cones and in particular are not applicable to $C^*$-algebras that contain nonzero projections.

A C*-algebra $A$ is called
 \emph{residually finite-dimensional} (abbreviated \emph{RFD}) if the finite dimensional
 representations of $A$ separate the points of $A$.
 We have shown in \cite{Dad:kk-top} that a separable, exact, {RFD} C*-algebra which satisfies the UCT is AF-embeddable and in fact it even embeds in a UHF-algebra $\bigotimes_{n=1}^\infty M_{k(n)}(\mathcal{C})$.
In the present note we point out that the arguments of \cite{Dad:kk-top} can be adapted to obtain a UHF-embeddability result for RFD $C^*$-algebras which does not assume the UCT but requires (local) finite generation of the even K-homology, see Definition~\ref{def-k-homology}.
\begin{theorem}\label{main}
  Let $A$ be a separable exact residually finite dimensional C*-algebra. If the rational $K^0$-homology  of $A$ is locally finitely generated, then $A$ embeds in a UHF-algebra.
\end{theorem}
\section{Preliminaries}
Let $A$ be a  $C^*$-algebra. A family $\mathcal{D}$ of $C^*$-subalgebras of $A$ is called exhaustive
if  for any finite subset $\mathcal{F}\subset A$ and any $\ep>0$ there exists $D \in \mathcal{D}$ such that $\mathcal{F}\subset_{\ep} D,$ i.e.  for each $x\in \mathcal{F}$ there is $d\in D$ such that $\|x-d\|<\ep$.
\begin{definition}\label{def-k-homology} We say that the  $K^0$-homology of $A$ is \emph{locally finitely generated} if
there is an exhaustive family $\mathcal{D}$ of $C^*$-subalgebras of $A$ such that for every $D \in \mathcal{D}$
the abelian group $K^0(D)=KK(D,\mathbb{C})$ is finitely generated. If instead we require the weaker condition that
each $\Q$-vector space $K^0(D)\otimes_\Z \mathbb{Q}$ is finite dimensional, then we say that the
 \emph{rational $K^0$-homology} of $A$ is {locally finitely generated}. The case when the vector space $K^0(A)\otimes_\Z \mathbb{Q}$ is itself finite dimensional is an obvious first example.
\end{definition}

 Let  $\mathcal{L}(\mathcal{H})$ denote the  linear operators
acting on a separable Hilbert space $\mathcal{H}$ and let $\mathcal{K}(\mathcal{H})$
denote the compact operators. We make the identifications
$\mathcal{L}(\mathbb{C}^k)=\mathcal{K}(\mathbb{C}^k) \cong M_k(\mathbb{C})$.
The unitary group of $M_k(\mathbb{C})$ is denoted by $U(k)$.
 Let
$A$ be a  $C^*$-algebra,  let $\mathcal{F} \subset A$
be a finite subset and  let $\ep >0$. If $\varphi:A \to \mathcal{L}(\mathcal{H}_{\varphi})$ and $\psi:A
\to \mathcal{L}(\mathcal{H}_{\psi})$
 are two maps, we write $\varphi\sim_{\mathcal{F},\ep} \psi$ if
there is a unitary $v:\mathcal{H}_\varphi \to \mathcal{H}_\psi $ such that
$\|v\varphi(a)v^*-\psi(a)\|< \ep$ for all $a \in \mathcal{F}$.

If $m$ is a positive integer and $\pi$ is a representation, then
 $m\pi$ will denote the representation $\pi \oplus \cdots \oplus \pi$ (m-times).
 The infinite direct sum  $\pi \oplus \pi \oplus \cdots$ is denoted by $\pi_\infty$.

We need the following approximation result.
 \begin{proposition}[{\cite[Prop.6.1]{Dadarlat:MRL}}]\label{prop:admissible-exists}
Let $A$ be a unital separable exact $C^*$-algebra and let $(\chi_n)_{n\geq 1}$ be a sequence
of unital representations of $A$ that separates the elements of $A$ and such that each representation in the sequence repeats itself
infinitely many times. For any
$\mathcal{F}\subset A$ a finite subset and any $\varepsilon >0$ there is an integer $r\geq 1$ such that if
$\pi=\chi_{1}\oplus \chi_{2}\oplus \cdots \oplus \chi_{r}$, then
for   any unital faithful representation $\sigma:A \to \mathcal{L}(\mathcal H)$
 with $\sigma(A) \cap \mathcal{K}(\mathcal H)=\{0\}$, one has $\sigma \sim_{\mathcal{F},\ep} \pi_\infty$.
\end{proposition}
\begin{definition}\label{admissible}
Let $A$ be a unital  RFD $C^*$-algebra. Let $\mathcal{F} \subset A$ be a
finite subset and let $\ep >0$. A unital representation $\pi:A \to
M_k(\mathbb{C})$ is called $(\mathcal{F},\ep)$-{admissible} if there is a
unital faithful representation $\sigma:A \to \mathcal{L}(\mathcal H)$ with $\sigma(A) \cap
\mathcal{K}(\mathcal H)=\{0\}$ ($\mathcal{H}=\mathbb{C}^k \oplus \mathbb{C}^k \oplus \cdots $)
such that
\begin{equation} \label{ad}\|\sigma(a)-\pi_\infty(a)\|<\ep, \quad \text{for all}\,
 a \in \mathcal{F}.\end{equation}
\end{definition}
\begin{remark}\label{rem}
  Note that if $\pi$ is $(\mathcal{F},\ep)$-admissible, then so is
$\pi \oplus \alpha$ for any unital finite dimensional
representation $\alpha$. Moreover $\|\pi(a)\| > \|a\|-\ep$ for
$a \in \mathcal{F}$. If a unital $C^*$-algebra $A$ is  separable exact and RFD,
then  Proposition \ref{prop:admissible-exists} guaranties the existence of $(\mathcal{F},
\ep)$-admissible representations for any finite set $\mathcal{F} \subset A$ and
any $\ep>0$.
\end{remark}

The following proposition is crucial for our embedding result.
It is based on a uniqueness theorem from \cite{DadEil:AKK}.
\begin{proposition}\label{prop:uniqueness}
Let $A$ be a unital separable exact RFD $C^*$-algebra. Let $\mathcal{F}
\subset A$ be a finite subset  and let  $\ep > 0$. Then for any
$(\mathcal{F},\ep)$-admissible representation
 $\pi:A \to M_k(\mathbb{C})$ and any
 two  unital representations $\varphi, \psi :A \to M_r(\mathbb{C})$, such that $[\varphi]=[\psi]\in K^0(A)$,
 there exist a positive integer $M$ and a unitary $u
\in U(r+Mk)$ such that
\begin{equation}\label{estimate}
\|u(\varphi(a) \oplus M \pi(a)) u^*-\psi(a)\oplus M\pi(a)\| < 3\ep,
\quad \text{for all}\,\,  a \in \mathcal{F}.
\end{equation}
\end{proposition}
\emph{Proof}. Fix $\mathcal{F}$, $\ep$ and $\pi$.
  Let $\sigma$ be a unital faithful representation of $A$  given by Definition~\ref{admissible}.
  In particular, $\sigma$ satisfies \eqref{ad} and $\sigma(A) \cap
\mathcal{K}(\mathcal H)=\{0\}$.
  Since $[\varphi]=[\psi]\in K^0(A)$, it follows that if we set $\Phi=\varphi\oplus \sigma$ and $\Psi=\psi\oplus \sigma$, then $\Phi(a)-\Psi(a)$ is a compact operator for all $a$ and the class of the Kasparov triple $(\Phi,\Psi,1)$ in $KK(A,\mathbb{C})=K^0(A)$ vanishes since
  \[ [\Phi,\Psi,1]=[\varphi\oplus \sigma, \psi\oplus \sigma, 1_r \oplus 1_{\mathcal{H}}]=[\varphi]-[\psi]=0.\]
  Moreover, both $\Phi$ and $\Psi$ are faithful representations whose images do not contain nonzero compact operators.
  This enables us to apply \cite[Thm.~3.12]{DadEil:AKK} and obtain that $\Phi$ is asymptotically unitarily equivalent to $\Psi$ via a continuous
  path of unitaries which are compact perturbations of the identity. In particular,
  there is a unitary $v\in U(\mathbb{C}^r\oplus \mathcal{H})$ of the form $v=1+x$ with $x\in\mathcal{K}(\mathbb{C}^r\oplus \mathcal{H})$ such that
 $\|v(\varphi(a) \oplus \sigma(a)) v^*-\psi(a)\oplus \sigma(a)\| < \ep$, for all
$a \in \mathcal{F}$. Using \eqref{ad} we obtain that
\begin{equation}\label{eqn:intermmediary}
  \|v(\varphi(a) \oplus \pi_\infty(a)) v^*-\psi(a)\oplus \pi_\infty(a)\| < 3\ep, \quad \text{for all}\,\,  a \in \mathcal{F}.
\end{equation}
 Since $\pi$ is a unital representation, it follows that the sequence of  projections $p_n=1_r\oplus n\pi(1)$ forms an approximate unit of $\mathcal{K}(\mathbb{C}^r \oplus \mathcal{H})$ and hence $[p_n,v]=[p_n,x]\to 0$ as $n\to \infty$. Since  each $p_n$ commutes with both $\varphi(a) \oplus \pi_\infty(a)$ and
$\psi(a) \oplus \pi_\infty(a)$,  we obtain from \eqref{eqn:intermmediary} that
\[\|(p_nvp_n)(\varphi(a) \oplus n\pi(a)) (p_nvp_n)^*-\psi(a)\oplus n\pi(a)\| < 3\ep, \]
for all $a \in \mathcal{F}$ and all sufficiently large $n$.  Moreover
 one can perturb the almost unitary operator $p_nvp_n$ to a unitary $u$ satisfying \eqref{estimate} for a sufficiently large value of $n$ denoted $M$.

\section{Proof of Theorem~\ref{main}}
Without any loss of generality, we may assume that $A$ is unital. We denote by $\mathrm{Rep_{fd}}(A)$ the set of unital finite dimensional representations of $A$.
  Since $A$ is separable and RFD, there is a sequence $(\chi_n)_{n\geq 1}$ in $\mathrm{Rep_{fd}}(A)$ which separates the points of $A$ and such that each representation in the sequence repeats itself
infinitely many times.
   Let $(x_n)_{n=1}^\infty$ be a dense sequence of elements of $A$  and let $\ep_n=2^{-n}$. Since $K^0(A)\otimes_\Z \Q$ is locally finitely generated, for each $n\geq 1$ there is a unital $C^*$-subalgebra $A_n$ of $A$ such that the $\Q$-vector space $K^0(A_n)\otimes_\Z \Q$ is finite dimensional and $X_n:=\{x_1,...,x_n\}\subset_{\ep_n} A_n$.
Fix a finite set $\mathcal{F}_n=\{a_{n,1},...,a_{n,n}\}\subset A_n$ such that $\|x_i-a_{n,i}\|<\ep_n$ for all $1\leq i \leq n$. Since $K^0(A_n)\otimes_\Z \Q$ is finite dimensional, its subspace $V_n$ generated by all the classes $\{[\chi_i|_{A_n}]\otimes 1: i\geq 1\}$ must also be finite dimensional.
 Thus there is an integer $r_n$ such that $V_n$ is generated by just $\{[\chi_i|_{A_n}]\otimes 1: 1\leq i \leq r_n\}$.

Define $\pi_n \in \mathrm{Rep_{fd}}(A)$ by $\pi_n=\chi_1\oplus \chi_2\oplus...\oplus \chi_{r_n}$.
By Proposition~\ref{prop:admissible-exists}, after increasing $r_n$, if necessary, we can moreover arrange that $\pi_n|_{A_n}$ is $(\mathcal{F}_n,\ep_n)$-admissible.

With these choices, we are going to construct a sequence of unital representations $\gamma_n :A \to M_{k_n}(\mathbb{C})$
 such that for all $n\geq 1$:
\begin{enumerate}[(i)\,\,]
  \item $k_{n+1}=k_nm_n$ for some positive integer $m_n$,
  \item $\gamma_n$ is unitarily equivalent to $\pi_n \oplus \alpha_n$ for some $\alpha_n \in \mathrm{Rep_{fd}}(A)$,
  \item $\|\gamma_{n+1}(x)- m_n\gamma_n(x)\|<5\ep_n$, for all $x\in X_n$.
\end{enumerate}
We will see that in fact each $\gamma_n$ is unitarily equivalent to a representation of the form  $q_1\chi_1\oplus q_2\chi_2\oplus...\oplus q_{r_n}\chi_{r_n}$, for integers $ q_i\geq 0$.

Set $\gamma_1=\pi_1\oplus \chi_1$. Suppose now that $\alpha_i$ and $\gamma_i$ were constructed for all $i\leq n$ such that the properties (i), (ii) and (iii) are satisfied. We construct $\alpha_{n+1}$ and  $\gamma_{n+1}$ as follows.

We need the following elementary observation. Suppose that $G$ is an abelian group such that the vector space $G\otimes_\Z \Q$ is finite dimensional and it is spanned by  $g_1\otimes 1,...,g_r\otimes 1$ with $g_i\in G$. Then for any $g\in G$ there are strictly positive integers $p,m, q_1,...,q_r,$ such that
\[pg+q_1g_1+\cdots +q_r g_r=m(g_1+\cdots  +g_r).\]

By applying this observation to the abelian subgroup of $K^0(A_n)$ generated by $\{[\chi_i|_{A_n}]: i\geq 1\}$
with $g_i=[\chi_i|_{A_n}]$, $i=1,...,r_n$, one obtains strictly positive integers $p,m, q_1,...,q_{r_n},$ such that
\[p[\pi_{n+1}|_{A_n}]+\sum_{i=1}^{r_n} q_i[\chi_i|_{A_n}]=m\sum_{i=1}^{r_n} [\chi_i|_{A_n}], \,\, \text{in}\,\, K^0(A_n).\]
Set $\alpha_{n+1}'=(p-1)\pi_{n+1}\oplus \big(\bigoplus_{i=1}^{r_n} q_i \chi_i\big) \oplus m \alpha_n$. Then
 \[[(\pi_{n+1}\oplus \alpha_{n+1}')|_{A_n}]=\sum_{i=1}^{r_n} m [\chi_i|_{A_n}]+m[\alpha_n|_{A_n}]=m[\pi_n|_{A_n}]\oplus m[\alpha_n|_{A_n}],\]
and hence
\([(\pi_{n+1}\oplus \alpha_{n+1}')|_{A_n}]=m [\gamma_n|_{A_n}],\) using (ii).
In particular, the representations $\pi_{n+1}\oplus \alpha_{n+1}'$ and $m\gamma_n$ have the same dimension.
Since $\pi_n|_{A_n}$ is $(\mathcal{F}_n,\ep_n)$-admissible, so is $\gamma_n|_{A_n}$ as noted in Remark~\ref{rem}. By Proposition~\ref{prop:uniqueness} applied to $A_n$, there is an integer $M\geq 1$ such that
\begin{equation*}
 (\pi_{n+1}\oplus \alpha_{n+1}')|_{A_n}\oplus M \gamma_n|_{A_n} \sim_{\mathcal{F}_n,3\ep_n} m \gamma_n|_{A_n} \oplus M \gamma_n|_{A_n}.
\end{equation*}
Set $\alpha_{n+1}=\alpha'_{n+1}\oplus M\gamma_n$, $m_n=m+M$ and $\gamma_{n+1}=\pi_{n+1}\oplus \alpha_{n+1}$.
Then
\begin{equation}\label{eqn:key}
 \gamma_{n+1}|_{A_n} \sim_{\mathcal{F}_n,3\ep_n} m_n \gamma_n|_{A_n}.
\end{equation}

Since  $\|x_i-a_{n,i}\|<\ep_n$ for all $1\leq i \leq n$, we deduce immediately from \eqref{eqn:key} that $\gamma_{n+1}\sim_{X_n, 5\ep_n} m_n\gamma_n$.
By conjugating $\gamma_{n+1}$ by a suitable unitary, we can arrange that $\|\gamma_{n+1}(x)- m_n\gamma_n(x)\|<5\ep_n$, for all $x\in X_n$.

Consider the UHF algebra $B=\underset{\longrightarrow}{\lim}\,\,M_{k(n)}(\mathbb{C})$ and let $\iota_n: M_{k(n)}(\mathbb{C}) \to B$ be the canonical inclusion.  Having the sequence $\gamma_n$
available, we construct a unital embedding $\gamma:A \to B$
 by defining $\gamma(x)$, $x \in \{x_1,x_2,\dots \}$, to be
the limit of the Cauchy sequence $(\iota_n\gamma_n(x))_{n\geq 1}$ and then
extend to $A$ by continuity.
 Note that $\gamma$ is a $*$-homomorphism, since all $\gamma_n$ are $*$-homomorphisms. Moreover, $\|\gamma(x)\|=\|x\|$ for all $x\in A$ since
$\|\gamma_n(a)\|> \|a\|-\ep_n$ for $a \in \mathcal{F}_n$ (by Remark~\ref{rem}) hence
$\|\gamma_n(x_i)\|\geq \|x_i\|-3\ep_n$,  as
$\|a_{n,i}-x_i\| <\ep_n$ for $1 \leq i \leq n$.
\begin{remark} It is clear from the proof that the conclusion of Theorem~\ref{main} holds under the weaker assumption that $A$ admits a separating sequence of finite dimensional representations $(\chi_n)_{n=1}^\infty$ such that for some exhaustive family $\mathcal{D}$ of $C^*$-subalgebras of $A$, the vector subspace of $K^0(D)\otimes_{\Z} \Q$
spanned by $\{[\chi_n|_D]\otimes 1\colon n\geq 1\}$ is finite dimensional for all $D\in \mathcal{D}$.
For example, this condition is satisfied if $A$ is the  suspension of a separable exact RFD $C^*$-algebra. Indeed, in that case one can choose a separating sequence $(\chi_n)_{n=1}^\infty$ with the property that $[\chi_n]=0$ in $K^0(A)$ for all $n\geq 1$.
\end{remark}
\textbf{Acknowledgements}.
The author would like to thank the referee  for a close
reading of the paper and for useful
suggestions.


\begin{thebibliography}{10}

\bibitem{Brown-Ozawa}
N.~P. Brown and N.~Ozawa.
\newblock {\em {$C^*$}-algebras and finite-dimensional approximations},
  volume~88 of {\em Graduate Studies in Mathematics}.
\newblock American Mathematical Society, Providence, RI, 2008.

\bibitem{Dadarlat:MRL}
M.~Dadarlat.
\newblock Residually finite dimensional {$C^*$}-algebras and subquotients of
  the {CAR} algebra.
\newblock {\em Math. Res. Lett.}, 8(4):545--555, 2001.

\bibitem{Dad:kk-top}
M.~Dadarlat.
\newblock On the topology of the {K}asparov groups and its applications.
\newblock {\em J. Funct. Anal.}, 228(2):394--418, 2005.

\bibitem{DadEil:AKK}
M.~Dadarlat and S.~Eilers.
 Asymptotic unitary equivalence in {$KK$}-theory.
{\em $K$-Theory}, 23(4):305--322, 2001.

\bibitem{Jamie-Gabe:connective}
J. Gabe.
\newblock {Traceless AF embeddings and unsuspended $E$-theory}.
 { arXiv:math/1804.08095} [math.OA], 2018.

\bibitem{Kir-icm}
E.~Kirchberg.
\newblock {Exact ${\rm {C}}\sp *$-algebras, tensor products, and the
  classification of purely infinite algebras}.
\newblock In {\em Proceedings of the International Congress of Mathematicians,
  Vol.\ 1, 2 (Z\"urich, 1994)}, pages 943--954, Basel, 1995. Birkh\"auser.

\bibitem{Ozawa:qd}
N.~Ozawa.
\newblock Homotopy invariance of {AF}-embeddability.
\newblock {\em Geom. Funct. Anal.}, 13(1):216--222, 2003.

\bibitem{ORS:elem-amen}
N.~Ozawa, M.~R{\o}rdam, and Y.~Sato.
\newblock Elementary amenable groups are quasidiagonal.
\newblock {\em Geom. Funct. Anal.}, 25(1):307--316, 2015.

\bibitem{Rordam:AF}
M.~R{\o}rdam.
\newblock A purely infinite {AH}-algebra and an application to
  {AF}-embeddability.
\newblock {\em Israel J. Math.}, 141:61--82, 2004.

\bibitem{RosSho:UCT}
J.~Rosenberg and C.~Schochet.
\newblock {The K\"unneth theorem and the universal coefficient theorem for
  Kasparov's generalized $K$-functor}.
\newblock {\em Duke Math. J.}, 55(2):431--474, 1987.

\bibitem{TWW:quasidiagonality}
A.~Tikuisis, S.~White, and W.~Winter.
\newblock Quasidiagonality of nuclear {$C^\ast$}-algebras.
\newblock {\em Ann. of Math. (2)}, 185(1):229--284, 2017.

\bibitem{Voi:qd}
D.~Voiculescu.
\newblock {A note on quasidiagonal $C^*$-algebras and homotopy}.
\newblock {\em Duke Math. J.}, 62(2):267--271, 1991.

\end{thebibliography}
\end{document}